\documentclass[a4paper,10pt]{article}
\def \Z {\mathbf Z}

\def \N {{\mathbf {N}}}
\def \Z {{\mathbf {Z}}}

\def \P {\mathbf {P}}
\def \B {{\cal B}}

\def\uu{\bigsqcup}
\def\eps{\varepsilon}
\def\CC{{\bf  C}}

\title{\bf Самоприсоединения  и типичные расширения эргодических систем}
\author{В.В.Рыжиков}
\date{}
\usepackage[T1]{fontenc} % иногда T1 вместо T2A
\usepackage[cp1251]{inputenc}  
\usepackage[russian]{babel}
\usepackage{graphicx}
\graphicspath{{pictures/}}
\DeclareGraphicsExtensions{.pdf,.png,.jpg}
\begin{document}
\large
\maketitle
\begin{abstract}  
В статье доказано, что типичные расширения динамической системы наследуют 
тривиальность самоприсоединений с попарной независимостью. Это  свойство связано с некоторыми известными  задачами
теории джойнингов и знаменитой проблемой Рохлина о кратном перемешивании.

Библиография: 15 названий, УДК 517.987  

\vspace{2mm}
Ключевые слова и фразы: \it самоприсоединения с попарной независимостью, 
локальная жесткость коциклов, типичные расширения действий, относительное кратное перемешивание.

\end{abstract} 

\section{Введение}  

В статье  рассматриваются сохраняющие вероятностную меру обратимые преобразования (автоморфизмы)
 пространства $(X,\B, \mu)$, образующие группу $Aut(\mu)$,  и их расширения, действующие в пространстве   $(X\times Y,\B\otimes\B, \mu\otimes\mu)$ в случае   $Y=X$.  
 Через $Ext(S)$ обозначены все косые произведения $R$  над $S$, иначе говоря, расширения преобразования $S$. Напомним, что косое произведение $R$ определено  формулой
$$R(x,y)= (Sx, R_x y), \ x\in X, \ y\in Y,$$
где $\{R_x\}$ -- измеримое семейство автоморфизмов пространства $(Y,\B,\mu)$. 
Метрика Халмоша на $Aut(\lambda)$,  $\lambda=\mu\otimes\mu$, задана формулой
$$ \rho(P,R)=\sum_i 2^{-i}\left(\lambda(PA_i\Delta RA_i)+\lambda(P^{-1}A_i\Delta R^{-1}A_i)\right),$$
где $\{A_i\}$  -- некоторое  фиксированное семейство,  плотное в алгебре $\B\otimes \B$. 
Оснащая    пространство  $Ext(S)$ метрикой $\rho$, которая полна, говорим, что
свойство  расширения \it типично, \rm если все расширения с этим свойством  содержат $G_\delta$-множество, плотное  в $Ext(S)$. Свойство автоморфизма $S$ \it стабильно, \rm если оно типично в $Ext(S)$.

Известно, что  бернуллиевость и $K$-свойство стабильны \cite{GTW}, также  установлена  стабильность  сингулярности спектра преобразования и свойства перемешивания \cite{22}. 
В \cite{AGTW} доказано, что  для эргодического автоморфизма $S$ с положительной
энтропией типичные расширения не только сохраняют все его инвариантны, но и оказываются изоморфными самому $S$.  В связи с этим   результатом  возникает ряд   вопросов. Если  эргодический автоморфизм с нулевой энтропией, например,   имеет корни, верно ли, что  его типичное расширение наследует это свойство? Если преобразование не имеет нетривиальных собственных инвариантных сигма-подалгебр (факторов),
верно ли, что его типичные расширения обладают лишь одним собственным  фактором? (Тривиальные факторы $\{X,\phi\}$, $\{X\times Y,\phi\}$  не учитываются.)  Более общий вопрос:  пусть  эргодическое действие некоторой группы    имеет в точности $n$ факторов, включая само действие, можно ли сказать, что  его типичное расширение обладает  $2n$ факторами? 

 Цель предлагаемой статьи -- добавить к списку стабильных свойств  инвариант, рассмотренный  дель Джунко и Рудольфом  в  работе \cite{JR}.  Дадим его определение.
Попарно независимым самоприсоединением порядка $n>2$ действия $\Psi$ называется мера на кубе $X^n$, для которой  проекции  на все двумерные грани  равны $\mu\otimes\mu$,  
причем мера инвариантна относительно диагонального действия 
 произведения $\Psi\times\dots\times \Psi$ ($n$ сомножителей). 

JR-свойство действия означает, что   для всех $n>2$ мера $\mu^n$ является  единственным 
его попарно независимым самоприсоединением   порядка $n$. Оно сохраняется при изоморфизме. Сформулируем основной результат статьи.
 
\vspace{2mm}
 \bf Теорема 1.1. \it  Типичные расширения наследуют JR-cвойство автоморфизма.\rm 

\vspace{2mm}
  Отсутствие нетривиальных самоприсоединений с попарной независимостью
 легко устанавливается для групповых действий, в слабом замыкании которых есть оператор $\eps I+(1-\eps)\Theta$, $\eps\in (0,1)$, где $I$ -- тождественный оператор, $\Theta$ --ортопроекция на пространство констант (см. \cite{R}).   В ряде случаев доказательство JR-свойства требует значителных усилий. Например, в \cite{FK} для класса потоков на поверхностях установлен  аналог  свойства Ратнер, этот инвариант влечет за собой так называемую квазипростоту порядка 2, которая в случае потоков запрещает нетривиальные самоприсоединения с попарной независимостью \cite{RT}. Как мы покажем, типичные расширения  преобразований, входящих в упомянутые потоки,  наследуют JR-свойство, а оно   влечет за собой  свойство кратного перемешивания.  Наш основной инструмент -- локальная жесткость расширений, применявшаяся в работе \cite{R97}. Это  свойство обычно  называют рекуррентностью коциклов (см.,
например, \cite {LP},\cite{L}). Такое  название не полностью отражает  специфику наших коциклов, для которых  рекуррентность   наблюдается на множестве меры, стремящейся к нулю, а на остальном множестве преобладает относительное перемешивание.  Мы  говорим, что соответствующие  коциклы обладают  $(I,\Theta)$-свойством. Здесь  тождественный оператор  $I$ символизирует локальную жесткость, а оператор 
ортопрекции на константы $\Theta$ -- относительное перемешивание.
%%%%%%%%%%%%%%%%%%%%%%%%%%%%%%%%%%%%%%%%%%%%%%%%%%%%%%%%%%%%%%%%%%%%%%%%%%%%%%%%%
\section{Свойство дель Джунко-Рудольфа } 
 Джойнингом набора действий $\Psi_1,\dots,\Psi_n$ называется мера на $X^n =X_1\times \dots\times  X_n$  ($X_i=X$), проекции которой на ребра  куба $X^n$ равны $\mu$, причем мера инвариантна относительно 
диагонального действия произведения $\Psi_1\times \dots\times  \Psi_n$. Если $\Psi_1,\dots,\Psi_n$ суть копии одного
действия, такой джойнинг называется самоприсоединением.

Говорим, что действие $\Psi$  принадлежит классу $S(m, n)$, $n>m>1$  (или обладает свойством $S(m,n)$), если всякое cамоприсоединение  порядка $n>2$ такое, что все проекции на  $m$-мерные грани  куба $X^n$ равны $\mu^m$,  является тривиальным,  т.е. совпадает с  $\mu^n$. 

Класс действий с  JR-свойcтвом эквивалентным образом определяется как  
$$S(2,\infty):=\bigcap_{n>2}S(2,n)=\bigcap_{n>2}S(n-1,n).$$ 
Авторы \cite{JR} доказали, что JR-свойcтво замкнуто относительно прямых произведений и  сохраняется при слабо перемешивающих компактных  расширениях автоморфизмов. Также они отметили, что 
для перемешивающих коммутативных действий JR-свойство влечет за собой \rm кратное перемешивание. \rm
Напомним, что автоморфизм $S$ перемешивает с кратностью $n$, если   для любых $A_0,A_1,\dots,A_n\in\B$ при $k_1,\dots, k_n\to\infty$ выполнено
$$\mu\left(A_0\cap S^{k_1}A_1\cap S^{k_1+k_2}A_2\cap \dots\cap  S^{k_1+\dots+k_n}A_n\right)\to 
\mu(A_0 )\mu(A_1 )\dots \mu(A_n).$$
%Это понятие введено Рохлиным в \cite{Ro}.

Кинг \cite{K} доказал, что $S(2,4)=S(2,\infty).$   В \cite{R} его замечательный   результат был дополнен автором:
$$S(2,\infty)=S(2m-1,2m), \ m>1.$$ Упомянутые утверждения универсальны в том смысле, что  верны для действий произвольных бесконечных групп. 

Упомянутые свойства  можно переформулировать в терминах факторов.
Фактором называется ограничение действия на инвариантную сигма-подалгебру.
JR-свойство действия $\Psi$ означает, что  всякий набор его попарно независимых факторов, 
изоморфных $\Psi$, является глобально независимым (действие, порожденное таким набором факторов, изоморфно прямому произведению этих факторов).  
Кинг в \cite{K} показал, что попарно независимый джойнинг действия  $\Psi$ со свойством $S(2,4)$
и двух копий произвольного действия $\Phi$ влечет за собой  тривиальность джойнинга. 
 В \cite{R93}, \S 3, этот результат  был усилен.

 \vspace{2mm}
\bf Теорема 2.1. \it Если фактор $\Psi$ обладает JR-свойством, попарная независимость  
факторов $\Psi$,$\Phi$,$\Pi$ влечет за собой  их глобальную независимость. \rm

\vspace{2mm}
\bf Следствие. \it Пусть фактор $\Psi$ обладает JR-свойством и независим с каждым из факторов $\Psi_1, \dots,\Psi_n$.
Если последние глобально независимы, то они глобально независимы вместе с фактором $\Psi$. \rm

\vspace{2mm}
Для доказательства следствия  сперва положим $\Phi=\Psi_1$, $\Pi=\Psi_2$. Из теоремы вытекает независимость $\Psi$ и 
$\Psi_1\bigvee \Psi_2$ (минимальный фактор, содержащий факторы $\Psi_1$ и
 $\Psi_2$).  Далее   рассматриваем $\Phi=\Psi_1\bigvee \Psi_2$, $\Pi=\Psi_3$ и получаем глобальную независимость факторов  $\Psi,$ $\Psi_1,\Psi_2$, $\Psi_3$. 
 Применяя теорему последовательно к факторам   
$\Phi=\Psi_1\bigvee$ $\dots$ $\bigvee \Psi_{k-1}$  и  $\Pi=\Psi_k$, 
 получаем утверждение следствия.

В теории джойнингов групповых действий  остался  открытым вопрос о совпадении  нечетных свойств $S(2m,2m+1)$ (напомним, что все  четные 
$S(2m-1,2m)$ при $m>1$ совпадают). Для некоторых некоммутативных действий нечетные свойства существенно слабее   четных.
Примерами служат  действия, порожденные  
автоморфизмами коммутативной группы и  сдвигами на этой группе. 

\vspace{2mm}
\bf  Действие с нечетным JR-свойством. \rm
Автоморфизмы группы $X= \Z_2\times\Z_2\times\Z_2\dots$  и групповые 
сдвиги на $X$ сохраняют меру Хаара $\mu$
на $X$.
Для  перестановки $\sigma$ натурального ряда $\N$ определим автоморфизм
$T_\sigma$
группы $X$ равенством $T_{\sigma}(\{x_i\})=\{x_{\sigma(i)}\}$.
Последовательности $\alpha=\{a_i\}\in X$ сопоставим
 сдвиг: $S_{\alpha}(x)=x+\alpha$, где $+$ обозначает  операцию сложения 
в $X$.   В качестве $\Psi$ рассмотрим действие, порожденное
всевозможными $T_\sigma$ и $S_{\alpha}(x)$. 

\vspace{2mm}
\bf Теорема 2.2. \it   Для всех $m> 0$  действие   $\Psi$
 обладает свойствами $S(2m, 2m+1)$, но оно не обладает JR-свойством.\rm

\vspace{2mm}
 Мера $\eta$, определенная формулой
$$
\eta(Y_1\times Y_2\times Y_3\times Y) =
\mu^3(\{(a,b,c): a\in Y_1, b\in Y_2, c\in Y_3,\,\, a+b+c\in Y\}),
$$
отлична от $\mu^4$,  инвариантна относительно диагонального действия 
$\Psi\otimes \Psi\otimes \Psi\otimes \Psi$,
а  ее проекции на 3-грани суть $\mu^3$. Поэтому $\Psi\notin S(3,4)$.
Несложное доказательство свойств $S(2m, 2m+1)$ дано в \cite{R96}, оно  использует тензорные произведения характеров группы $X$.  Действие  $\Psi$ содержит 
конечно-порожденные поддействия с аналогичным свойством.

\vspace{2mm}
\bf Вопрос дель Джунко-Рудольфа и проблема Рохлина. \it Существует ли  слабо перемешивающее $\Z$-действие с нулевой энтропией, которое не принадлежит классу $S(2,\infty)$? \rm Отрицательный  ответ на этот вопрос в классе перемешивающих преобразований одновременно будет  решением проблемы Рохлина: \it обладает ли  перемешивающий автоморфизм  свойством  перемешивания кратности $k>1$? \rm 

В \cite{JR} было отмечено, что для $\Z^2$-действий 
примерами перемешивающих систем без JD-свойства  служат известные перемешивающие $\Z^2$-действия  Ледрапье, не обладающие кратным перемешиванием (новый взгляд на такие действия изложен в   \cite{T}).

\section{ $(I,\Theta)$-коциклы сохраняют JR-свойство} 
Пусть мера $\nu$ на $Y\times Y_1\times\dots Y_n$  пректируется  на сомножитель 
$Y$ в меру $\mu$, а ее проекция на  $Y^n=Y_1\times\dots Y_n$ 
совпадает с  $\mu^n$.  Класс таких мер на $Y^{n+1}$ обозначим через $M_{n+1}$. 

\vspace{2mm}
{\bf Лемма 3.1.} \it  Пусть 
$R=(S,R_x)$ является  слабо перемешивающим расширением преобразования $S$ с $JR$-свойством.  Если расширение $R$ не обладает $JR$-свойством,  найдется  $\mu^n$-измеримая функция
 ${\cal M} : X^{4} \to M_{4}$, удовлетворяющая  тождеству 
$$
   {\cal M} (S(x),S(x_1),S(x_2),S(x_3))\equiv
(R_{x}\otimes R_{x_1}\otimes R_{x_2}\otimes R_{x_3}){\cal M}(x,x_1,x_2,x_3), \eqno (1)
$$
причем ${\cal M}$ существенно отлична от $\mu^4$.
\rm

Доказательство. При отсутствии $JR$-свойства для преобразования $R:X\times Y\to X\times Y$ найдется  $\nu$ -- нетривиальное самоприсоединение  порядка 4 со стандартными проекциями на 3-грани в кубе $\bar X^4$, $\bar X=X\times Y$. 
Проекция $\nu$ на $X^4$ совпадает с $\mu^4$,
так как проекция инвариантна относительно  $S^{\otimes 4}$, но $S$ обладает $JR$-свойством.
В силу сказанного  мера $\nu$ имеет представление 
$$\nu ((A\times B)\times (A_1\times B_1)\times (A_2\times B_2)\times (A_3\times B_3)) =
\int_{A\times A_1\times A_2\times A_3} \nu_{(x,x_1,x_2,x_4)}(B\times B_1\times B_2\times B_3)\,d\mu^4, $$
где $\nu_{(x,x_1,x_2,x_4)}$ -- соответствующая система условных мер.

В силу следствия из теоремы 2.1  из  независимости фактора $S$, обладающего JR-свойством,  
с тремя  координатными $R$-факторами в произведении $R\times R\times R$  вытекает независимость $S$-фактора с этим произведением. 
 Следовательно,
$$\nu ((A\times Y)\times (A_1\times B_1)\times (A_2\times B_2)\times (A_3\times B_3)=\mu(A)\mu(A_1)\mu(A_2)\mu(A_3)\,\mu(B_1)\mu(B_2)\mu(B_3).$$
Это означает, что  почти все условные меры $\nu_{(x,x_1,x_2,x_3)}$ 
лежат в классе $M_{4}$. \rm Лемма доказана.

Ниже определяется свойство расширений, которое, как мы покажем,  вынуждает функции,  удовлетворяющие (1), совпадать с константой $\mu^4$.    В \S 4 будет установлено, что типичные расширения обладают таким свойством.

\bf Определение $(I,\Theta)$-свойства коцикла. \rm 
 Пусть для  косого произведения $R=(S,R_x)$ и некоторой последовательности  $p_j\to\infty$  выполняются два условия.
Первое  -- относительное перемешивание по мере вдоль последовательности $p_j$: для любого $\eps>0$
$$\mu\left(x:\, dist_w(C(x,R^{p_j}),\,\Theta)<\eps\right)\to 1, \ j\to\infty,$$
где $dist_w$ обозначает некоторую фиксированную метрику, задающую слабую  сходимость операторов в $L_2(\mu)$.  Потребуем также сходимость $S^{p_j}\to_w\Theta$.

 Второе условие   -- локальная жесткость коцикла: 
для всякого множества $A\subset X$ положительной меры для почти всех  $x\in  A$ найдется последовательность  $j'\to\infty$ (она зависит от $x$) такая, что выполнены
$S^{p_{j'}}x\in  A$ и
$$ \rho\left(C(x,R^{p_{j'}}), Id\right)\to 0,\ j'\to\infty,\ 
C(x,R^p):=R_{S^{p-1}x}\dots R_{Sx}R_x$$
Здесь и всюду ниже $\rho$ обозначает  метрику Халмоша на пространстве $Aut(\mu)$.

При выполнении  указанных  условий  говорим, что  косое произведение $R$  обладает $(I,\Theta)$-коциклом (или говорим, что коцикл обладает $(I,\Theta)$-свойством).

\vspace{2mm}
 \bf Теорема 3.2. \it  Расширение  $R=(S,R_x)$, обладающее  $(I,\Theta)$-коциклом,  наследует  JR-свойство автоморфизма $S$.  
\rm 

\vspace{2mm}

Напомним, что через $M_{4}$ в статье обозначен класс  мер на кубе $Y\times Y_1\times Y_2\times Y_3$,   проекции которых на ребро $Y$ и  на  грань
$Y_1\times Y_2\times Y_3$ равны $\mu$  и $\mu^3$, соответственно. 

\vspace{2mm}
 \bf  Лемма 3.3. \it  Пусть для последовательностей автоморфизмов $T_j,T_{1,j},T_{2,j},T_{3,j}$,
 имеют место сходимости
 $T_j\to I$ и  $T_{1,j},T_{2,j},T_{3,j}\to_w\Theta$ при $j\to\infty$.
Тогда для всякой меры $\nu\in M_4$ выполнено 
$$\nu_j=(T_j\times T_{1,j}\times T_{2,j}\times T_{3,j})\nu\to \mu^{4}.$$ \rm

\vspace{2mm}
Доказательство. Воспользуемся связью между 
мерами, являющимися полиморфизмами в терминологии \cite{V}, и  марковскими операторами. Отметим,
что слабую сходимость мер в нашем случае мы  отождествляем со слабой сходимостью соответствующих марковских операторов.
Связь меры с оператором  задается формулой 
$$(Pf\,,\, f_1\otimes f_2\otimes f_3)= 
\int_{Y^{4}} f\otimes f_1\otimes f_2\otimes f_3\, d\nu$$
(здесь подразумевается скалярное произведение в пространстве $L_2(\mu^3)$).
Из условий леммы имеем:  $PT_jf$ сходится к $Pf$ по норме, а  
последовательность 
$F_j=T_{1,j}f_1\otimes T_{2,j}f_2\otimes T_{3,j}f_3$ слабо сходится к константе $const\equiv c$, 
которая равна произведению интегралов функций $f_1,f_2,f_3$. 
Так как  для  оператора $P$ выполнено   $(Pf\,,\, const)=c\int_Y  \, f d\mu$ (напомним, что марковские операторы сохраняют интеграл от функции),
получим 
$$\int_{Y^{4}} f\otimes f_1\otimes f_2\otimes f_3\, d\nu_j \ =
(PT_jf\,,\, F_j)\,
\to\ \int_{Y^4} f\otimes f_1\otimes f_2\otimes f_3\, d\mu^4.$$
Лемма доказана.

Доказательство теоремы 3.2.
Пусть $dist$ --  расстояние между мерами класса $M_4$, совпадающее
с расстоянием $dist_w$ между  марковскими операторами, соответствующими этим мерам, 
а сходимость по  $dist_w$ эквивалентна слабой сходимости операторов.
Предположим, что   $\mu^4(U)>0$ для множества
$$ U=\{(x, x_1,x_2,x_3)\, :\, dist({\cal M}(x, x_1,x_2,x_3),\,\mu^4)> a\}$$
при некотором  $a>0$.  Обозначим
$$ U_x=\{(x_1,x_2,x_3)\, :\, (x, x_1,x_2,x_3)\in U\}.$$

Найдется множество $A\subset X$ положительной меры такое, что для некоторого $u>0$ и множества 
$V\subset X_1\times X_2\times X_3$ выполнены условия: для всех $x\in A$
$$\mu^3(U_x)>3u, \ \ \mu^3(U_x\Delta V)<u^2.$$

Фиксируем точку $x\in A$ с соответствующей последовательностью $p_{j'}$, которую будем
далее обозначать как $p_j$ (см.  определение $(I,\Theta)$- свойства коцикла). 
В силу того, что 
$$(S\otimes S\otimes S)^{p_j}\to_w\Theta\otimes\Theta\otimes\Theta,$$  для всех достаточно больших $j$ имеем
$$ \mu^3(V\bigcap (S\times S\times S)^{p_j} V)>4u^2.$$
Так как для всех $x\in A$ выполнено $\mu^3(U_x\Delta V)<u^2$,  при $S^{p_{j'}}x\in A$
получаем
$$ \mu^3(U_{S^{p_j}x}\bigcap (S\times S\times S)^{p_j}U_x)>2u^2.$$
Рассмотрим множества 
$$ V_{x,j}= U_x\bigcap (S\times S\times S)^{p_j}U_x\bigcap W_{j},$$
$$
W_j=\{(x_1,x_2,x_3)\,:\, dist_w(C(x_1,R^{p_j})\otimes C(x_2,R^{p_j}))\otimes 
C(x_3,R^{p_j})\,, \,\Theta\otimes\Theta\otimes\Theta)<{\eps_j}\},
$$
где  $dist_w$ -- расстояние между операторами на $L_2(\mu^3)$, задающее слабую операторную топологию.

Eсли  последовательность ${\eps_j}$ достаточно медленно стремится к 0, то в силу определения 
$(I,\Theta)$-свойства коцикла для нашего расширения выполняется   $ \mu^3(W_j)\to 1$.
Для достаточно больших $j$ имеем неравенства
$$ \mu^3(V_{x,j})>u^2, \ \ \mu^3(\bigcap_i\bigcup_{j>i}V_{x,j}) \geq u^2.$$
Mножество $\bigcap_i\bigcup_{j>i}V_{x,j}$ состоит из троек $(x_1,x_2,x_3)$, принадлежащих  множествам  $V_{x,j}$ для бесконечного набора $J$ индексов $j$. Фиксируем некоторую такую тройку. Для  четверки $(x,x_1,x_2,x_3)$  при условии, что $j\in J$ и $j\to\infty$, имеем 
$$C( x,R^{p_{j}})\to I, \ S^{p_{j}}x\in A,$$
$$C( x_1,R^{p_{j}}),\, C(x_2,R^{p_{j}}),\, C( x_3,R^{p_{j}})\to_w\Theta.$$

Из тождества $(1)$ следует равенство 
$${\cal M}(S^{p_j}x,S^{p_j}x_1,S^{p_j}x_2,S^{p_j}x_3)=
\CC(x,x_1,x_2,x_3,p_j){\cal M}(x,x_1,x_2,x_3),
$$
где  
$$\CC(x,x_1,x_2,x_3,p)= C(x,R^{p})\times C(x_1,R^{p})\times C(x_2,R^{p})\times C(x_3,R^{p}).$$

 Меры $\nu_j={\cal M} (S^{p_j}x,S^{p_j}x_1,S^{p_j}x_2,S^{p_j}x_3)$ при $j\in J$ лежат в $U,$
 следовательно, они отделены от $\mu^4$.
Но в силу  леммы 3.3 имеет место сходимость 
$$\nu_j=\CC(x,x_1,x_2,x_3,p_{j})\nu_0\to \mu^4,  \ j\in J, \, j\to\infty.$$
Предположение о том, что    $\mu^4(U)>0$, приводит к противоречию.
Таким образом, наличие $(I,\Theta)$-коцикла запрещает  решения (1), отличные от ${\cal M}(x,x_1,x_2,x_3)\equiv  \mu^4$. Теорема 3.2 доказана.

\section {Типичные расширения обладают $(I,\Theta)$-коциклом}
 
\bf Теорема 4.1. \it  Для слабо  перемешивающего  преобразования  $S$ в пространстве $Ext(S)$ найдется типичное множество $G$ такое, что для всякого $R=(S,R_x)$ из $G$ и всякого множества $A\subset X$ положительной меры    найдутся  последовательности $\delta_j\to 0$ и $p_j\to\infty$ такие, что $S^{p_j}\to_w\Theta$  и 
$$\mu\left(x':\, dist_w(C(x', R^{p_j}),\,\Theta)<\delta_j\right)\to 1,\eqno ( \ast) $$
где $ C(x, R^p)$ обозначает произведение $R_{S^{p-1}x}\dots R_{Sx}R_x$.

Для почти всех $x\in A$ найдутся  последовательности $j'\to\infty$ (они  зависят от $x$),   
для которых
  $$S^{p_{j'}}(x)\in A, \ \rho(C(x,R^{p_{j'}}),\, Id)\to 0,\eqno (\ast\ast)  $$ 
\rm

\vspace{2mm} 

Для доказательства того, что $(\ast\ast)$ выполняется для почти всех точек из $A$, достаточно убедиться в том, что всякое множество
 положительной меры  содержит хотя бы  одну точку, удовлетворяющую $(\ast\ast)$. Назовем такую точку хорошей.  Множество всех   точек из $A$, не являющихся хорошими,  имеет нулевую меру, так как иначе в нем найдется хорошая  точка.   

 Мы предъявим последовательность  $p_j\to \infty$ такую, что для нее выполнены неравенства
$$\mu\left(x':\, dist_w(C(x', R^{p_j}),\,\Theta)< 1/ j \right)> 1- 1/j, $$
и последовательность множеств $A_j\subset A$  положительной меры, для которой
 $A_{j+1}\cup S^{p_{j}}A_{j+1}\subset A_j$  и  для всех $x\in A_{j+1}$ выполнено
$$ \rho(C(x,R^{p_j}),\, Id)<  1/ j .  $$ 
Не ограничивая общности, будем считать, что $X$ является компактом $[0,1]$,
 а множества $A_j$ выбираются замкнутыми.  Их  пересечение будет содержать хорошую точку.  Таким образом, наша цель  -- для каждого элемента  $R$  некоторого фиксированного типичного подмножества $G\subset Ext(S)$  и каждого множества $A$, $\mu(A)>0$,  найти требуемые  последовательности $p_j$ и   $A_j$.  Построение  такого $G$  составляет  основную часть доказательства теоремы.

\bf  Выбор плотного семейства расширений. \rm   Все расширения автоморфизма $S$, когомологичные  произведению $S\times Id$, образуют плотный в $Ext(S)$ класс.  
 Они имеют вид  $J^{-1}(S\times Id)J$, где $J=(Id, J_x)$ для некоторого измеримого семейства автоморфизмов $\{J_x\}$.
В  классе $Ext(S)$ зафиксируем счетное плотное множество $\{R_i\}$,
 $$ R_i=J_i^{-1}(S\times Id)J_{i}, \ J_i= (Id, J_{i,x}),$$
где автоморфизмы  $J_{i,x}$ устроены следующим образом:
для каждого $i$  определено разбиение $\xi_i$ пространства $X$ на   непересекающиеся множества
$B_{i, k}$, $1\leq k\leq 2^i$, одинаковой меры, т.е. $\mu(B_{i,k})=2^{-i}$, 
причем  для  фикированого $k$ автоморфизмы $J_{i,x}$ одинаковы при   $x\in B_{i,k}$.
Разбиения $\xi_i$ выбираем таким образом, чтобы они стремились к разбиению на точки, причем для всех $i$ и всех $k\leq 2^i$ выполнялось   $B_{i,k} =B_{i+1,2k-1}\uu B_{i+1,2k}$.
Существование нужного   семейства $\{J_i\}$,  плотного в $Ext(Id)$,   вытекает из сепарабельности пространства $Aut(\mu)$ и стандартных  фактов  теории меры.

\bf  Выбор окрестностей элементов плотного множества. \rm Так как 
$$R_i^p(x,y)=(S^px,\,C(x,R_i^p) y)=(S^px,\, J^{-1}_{i,S^p(x)}J_{i,x}\, y),$$
при $x,S^px\in B_{i, k}$  имеем $J_{i,x}^{-1} J_{i,x'}= Id$, следовательно,
$$\rho(C(x,R_i^p),\, Id)=0.$$

Через $\hat B$ обозначим индикаторы множеств $B$,
 положим $$\P_L=\frac 1 {2L} \sum_{p=L+1}^{3L} S^p.$$
 В силу эргодичности преобразования $S$ при  $L\to\infty$ выполняется 
$$\left\|\P_L \hat B_{i,k} -\mu(B_{i,k}) \hat X\right \|_1\to 0$$
(следствие  классической теоремы о сходимости эргодических средних в $L_1(\mu)$).

Фиксируем последовательность $\eps_i$, полагая  $\eps_i=o(\mu(B_{i,k})^2)$. 
Для всякого $i$  найдется число  $L(i)$ такое, что для всех
$k$, $1\leq k\leq i$, выполнено   неравенство 
$$\left\|\P_{L(i)} \hat B_{i,k} -
\mu(B_{i,k}) \hat X\right \|_1 < \eps_i\mu(B_{i,k})^2, \eqno (4)$$
причем, с  учетом того, что $S$ слабо перемешивает,  для большинства  $p$, 
$L(i)<p\leq 3L(i)$,  будет обеспечено условие 
$$|\mu(B_{i,k}\cap S^pB_{i,k})-\mu(B_{i,k})^2|<\eps_i\mu(B_{i,k})^2.$$
Большинство означает в нашем случае, что его мощность не меньше, чем $(1-\eps_i)2L(i)$.

Для расширения  $R\in Ext(S)$ определим множества 
$$D(p, B,\delta, R)=\left\{x\in B\cap S^{-p}B\,:\, \rho(C(x,R^p),Id)<\delta\right\}.$$
   Непосредственно из определений при $\delta>0$ вытекает равенство
$$B_{i,k}\cap S^{-p}B_{i,k}=D(p,B_{i,k}, \delta,R_i).$$ 
%следовательно, для всех натуральных $L$
%$$\frac 1 L \sum_{p=1}^L \hat D(p,B_{i,k},1/j, R_i) - 
%\frac 1 L \sum_{p=1}^L \hat B_{i,k}\,S^p\hat B_{i,k}=0.      \eqno (5)$$

Выбираем такую окрестность $U(j,i)$ автоморфизма $R_i$, что для 
всякого ее элемента $R$ для всех $k$, $1\leq k\leq i$,  и всех $p$, 
$L(i)< p\leq 2L(i)$,    выполнены  неравенства 
$$\mu\left( D(p,B_{i,k},1/j, R_i)\, \Delta\, D(p,B_{i,k},1/j, R)\right)\,
 <\,\eps_i\mu(B_{i,k})^2. \eqno (5)$$

%Тогда с учетом  $(5)$ получим
%$$\left\|\frac 1 {2L(i)} \,\sum_{p=L(i)+1}^{3L(i)} \hat D (p,B_{i,k}, 1/j,R)-  \mu(B_%{i,k})\P_{L(i)}\hat B_{i,k}  \right\|_1<  \eps_i\mu(B_{i,k})^2.$$

\bf  Определение типичного множества. \rm Пусть  
$$ R\in\,  G={ \mathop{R\,WM}}\bigcap\bigcap_{j}\bigcup_{i>j}  U(j,i).$$
 Класс   ${\mathop{R\,WM}}\subset Ext(S)$ состоит из расширений, обладающих свойством относительного слабого перемешивания (соответствующие определения даны ниже).  Наша ближайшая цель -- показать, что для  всякого $A\subset X$, $\mu(A)>0$, и всякого $j$  найдется   такое $i$, что 
  $$ |\{p:\, L(i)<p\leq 3L(i), \ \ \mu(D(p,A, 1/j,R))>0\}|>L(i).\eqno (6)$$

Фиксируем положительное   $c<0.01$. Для некоторого $i(A)$  для всех $i> i(A)$  найдется   $k$, $1\leq k\leq i$, для которого выполнено     
$\mu(A\cap B_{i,k})>(1-c)\mu(B_{i,k}) $. 
Сказанное является следствием того, что диадические разбиения $\xi_i$ 
образуют монотонную последовательность и стремятся к разбиению на точки. Обозначим  
$$ B_0=B_{i,k}, \ \ A_0=A\cap B_0,\ \ \Delta_0=B_0\setminus A_0.$$

Поскольку  $$ A_0\cap S^p A_0\setminus D(p,A_0, 1/j,R)\,\subset \, 
B_0\cap S^p B_0\setminus D(p,B_0, 1/j,R)$$  имеем
$$\mu(A_0\cap S^p A_0\setminus D(p,A_0, 1/j,R))\,\leq \,
\mu(B_0\cap S^p B_0\setminus D (p,B_0,1/j, R))\, <\, \eps_i\mu(B_{i,k})^2. $$ 
Учитывая $(5)$,  получим
$$ \left\| \hat D(p,A_0, 1/j,R)  - \hat A_0\cap S^p \hat  A_0\right\|_1<  \eps_i\mu(B_{i,k})^2,$$
что влечет за собой 
$$\left\|\frac 1 {2L(i)} \,\sum_{p=L(i)+1}^{3L(i)} \hat D(p,A_0, 1/j,R)-
\hat A_0 \P_{L(i)} \hat A_0   \right\|_1<  \eps_i\mu(B_{i,k})^2.$$

 Из равенства  $\hat B_0=\hat A_0+ \hat\Delta_0$ вытекает
$$ \hat A_0 S^p \hat A_0 =
\hat B_0 S^p \hat B_0   -\hat \Delta_0 S^p \hat B_0 -
\hat A_0 S^p \hat \Delta_0.  $$
Из $(4)$ следует, что 
$$\left\| \hat \Delta_0 \P_{L(i)}  \hat B_0   \right\|_1< 2c \mu(B_0)^2,$$
$$\left\| \hat A_0 \P_{L(i)}  \hat \Delta_0    \right\|_1 = 
\left\|\hat \Delta_0  \hat \P_{L(i)} ^\ast A_0     \right\|_1 \leq
\left\|\hat \Delta_0   \P_{L(i)} ^\ast \hat   B_0    \right\|_1 =
\left\|\hat \Delta_0  T^{4L_i} \P_{L(i)}  \hat   B_0    \right\|_1<2c \mu(B_0)^2,  $$
где $\P_L^\ast=\frac 1 {2L} \sum_{p=L+1}^{3L} S^{-p}.$
Следовательно, 
$$\left\| \hat A_0 \P_{L(i)}  \hat A_0   \right\|_1=
\left\|   \hat B_0 \P_{L(i)}  \hat B_0   \right\|_1 
-\left\| \hat \Delta_0 \P_{L(i)}  \hat B_0   \right\|_1
-\left\|\hat A_0 \P_{L(i)}  \hat \Delta_0    \right\|_1 >(1-5c) \mu(B_0)^2.$$

Таким образом, 
$$\left\| \hat A_0 \P_{L(i)}  \hat A_0   \right\|_1\, >\, 0.9 \mu(B_0)^2,$$
следовательно, число  таких $p$, удовлетворяющих  условиям $L(i)<p\leq 3L(i)$ и 
 $\mu(A_0\cap S^pA_0)>c\mu(B_0)^2$, больше  $L(i)$.

Мы получили (6) для всякого $R$, принадлежащего  плотному $G_\delta$-множеству  $\bigcap_{j}\bigcup_{i>j}  U(j,i)$.
Теперь  рассмотрим типичное   свойство, которое обеспечит выполнение $(\ast)$ в утверждении теоремы.

\bf  Относительное слабое  перемешивание. \rm Класс   $R\,W M $ состоит из косых
произведений $(S,R_x)$ таких, что  $(S,R_x\times R_x)$ эргодично относительно меры 
$\mu^3$. Типичность $R\,W M $ в пространстве $Ext(S)$ доказана в  работе \cite{GW}.
Свойство относительного слабого перемешивания влечет за собой для  любых измеримых  $Y_1,Y_2\subset Y$  сходимость 
$$  \int_X \frac 1 L \sum_{p=1}^{3L}\left(\mu(C(x, R^p)Y_1\cap Y_2)-
\mu(Y_1)\mu(Y_2)\right)^2 d\mu(x)\to 0, \ \ L\to\infty. $$
  
Поэтому для всякого $\delta>0$  при  больших значениях $L$ верно, что для более, чем $(1-\delta)2L$  натуральных чисел $p$ из интервала   $(L,3L)$,  выполнено
$$\mu\left( x:  \, dist_w(C(x, R^p),\Theta)<\delta\right)>1-\delta.$$  

\bf Завершение доказательства теоремы 4.1. \rm Для всякого  расширения $R\in G$ мы установили следующее утвержение.  Пусть  $\delta_0>0$ и множество $A$ имеет  положительную  меру. Тогда  найдется  $A_0\subset A$, $\mu(A_0)>0$ и  целочисленный интервал $[L_1, 3L_1]$ ($L_1$ может быть выбрано сколь угодно большим), в котором  большинство элементов    $p$   удовлетворяют условию
$$\mu\left(x':\, dist_w(C(x', R^{p}),\,\Theta)<\delta_0\right)> 1-\delta_0, $$
и более половины элементов  $p$  удовлетворяют условию
$$ \mu (D(p,A_0, \delta_0, R))> c\mu(A_0)^2/2 ,$$
где, напомним,    
$$ D(p,A_0, \delta_0, R)=\left\{x\in A_0\,:\, S^{p}(x)\in A_0, \ \rho(C(x,R^{p}),\, Id)
< \delta_0\right\}.$$ 
 Выберем один из таких элементов $p_1$  и положим  $ A_1=D(p_1,A,\delta_0,R)$.
 
Фиксируем последовательность $\delta_i\to +0$
Теперь, применим предыдущее утверждение  к множеству  $A_1$ (вместо множества $A$), мы найдем множество $A_2\subset A_1$, потом найдем $A_3\subset A_2$   и т.д., в результате получим 
%$Aтогда для $\delta_1>0$ найдется  соответствующие  $p_2\in [L_2, 3L_2]$ и  множество $A_2=D%(p_1,A_1,\delta_1,R,\delta_1)$ положительной меры. Потом в качестве $A$ рассмотрим $A_2$,  последовательность $p_j$ и последовательность вложенных множеств $A_j$, для которых выполняется  
$$\mu\left(x':\, dist_w(C(x', R^{p_j}),\,\Theta)<\delta_j\right)> 1-\delta_j, $$
 и 
$$S^{p_j}(x)\in A_{j-1}, \ \rho(C(x,R^{p_j}),\, Id)< \delta_j $$ 
для всех $x\in A_j$. Из сказанного, как  пояснялось выше,  вытекают требуемые утверждения $(\ast),(\ast\ast)$. Теорема 4.1. и, следовательно, теорема 1.1 доказаны.

В   \S 2  обсуждалось  нечетное JR-свойство. Неизвестно, сохраняется ли оно 
при типичных расширениях. Отметим также, что  примеры автоморфизмов, не обладающих JR-свойством, но   обладающих нечетным JR-свойством, не обнаружены, хотя имеются  контрпримеры   среди действий некоторых некоммутативных групп.  

Отметим также, что из стабильности свойства перемешивания (см. \cite{22}) и теоремы 1.1
вытекает следующее утверждение.

 \bf Теорема 4.2. \it  Типичные расширения автоморфизма $S$ наследуют свойство кратного перемешивания, если $S$  обладает  JR-свойством.\rm 

\vspace{2mm}
\section{Кратное перемешивание и расширения} 
Сохраняют ли типичные расширения свойство кратного перемешивания?  Этот вопрос остался    открытым, так как для возникающего самоприсоединения
$\nu$ отсутствует обоснование того, что условные меры $\nu_{(x_0,x_1,\dots,x_n)}$ принадлежат классу $M_{n+1}$. Рассмотрим примеры, когда требуемая принадлежность устанавливается.

\vspace{2mm}
{\bf Теорема 5.1.} (\cite{R97})  \it Пусть $R:X\times Y\to X\times Y$ -- косое произведение над автоморфизмом $S$, заданное формулой
$$
  R(x,y)=(S(x),T^{n(x)}(y)), \quad \int n(x)d\mu = 0.
$$
Если $R,T$ обладают перемешиванием, а автоморфизм  $S$ перемешивает с кратностью $k$, то  косое произведение $R$ наследует  перемешивание  кратности $k$.\rm

\medskip
Особенность  косого произведения $R=(S, T^{n(x)})$ в том, что $R$  коммутирует с 
$Id\times T$, причем действие $T$ перемешивающее. В ходе доказательства, которое 
похоже на доказательство теоремы 3.2, эти свойства    обеспечивают  принадлежность
соответствующих условных мер классу $M_3$. Так как в \cite{R97}  соответствующая  аргументация о принадлежности условных мер к классам $M_n$  не была приведена,   восполним этот пробел сейчас.

Пусть $k=2$. Рассмотрим     меру $\nu$ вида
$$
\nu (C_0\times C_1 \times C_2) =
\lim_{j\to\infty}
 \lambda (C_0 \cap R^{z_1(j)}C_1  \cap R^{z_2(j)}C_2),\ \ z_1(j),z_2(j), |z_1(j)-z_2(j)|\to\infty.
$$
Нам нужно убедиться в том, что мера $\nu$ принадлежит  классу $M_3$.
Так как $R$ коммутирует с преобразованием $Q=Id\times T$, 
мера $\nu$ инвариантна относительно $Q\times Q\times Q$, ее
проекция на $Y\times Y \times Y$ инвариантна относительно $T\times T \times T$ 
и принадлежит классу $M_3$ (это вытекает из свойства перемешивания действия  $R$). 
Так как проекция меры $\nu$ на $X\times X \times X$ есть $\mu^3$, получим
представление
$$\nu ((A_1\times B_1)\times (A_2\times B_2)\times (A_3\times B_3)) =
\int_{A_1\times A_2\times A_3} \nu_{(x_1,x_2,x_3)}(B_1\times B_2\times B_3)\,d\mu^3, $$
где $\nu_{(x_1,x_2,x_3)}$ --  система условных мер на  на $Y\times Y \times Y$ .

Почти все условные меры   $\nu_{(x_1,x_2,x_3)}$ 
лежат в классе   $M_3$, так как их проекции на ребро $Y_1$ совпадают с $\mu$,
 а проекции на грань  $Y_2\times Y_3$  совпадают с эргодической мерой  $\mu\otimes\mu$  (следствие того, что эргодичность интеграла по инвариантным мерам означает, что почти все эти меры совпадают с этим интегралом). Нужное $(I,\Theta)$-свойство коцикла  следует из теоремы Аткинсона-Крыгина и свойства перемешивания автоморфизма $T$ (см. \cite{R97}), оно тривиализует систему условных мер, т.е.  $\nu_{(x_1,x_2,x_3)}\equiv \mu^3$.
Таким образом,  $\nu=\lambda^3$, и тем самым перемешивание кратности 2 установлено.  
Для  кратности $k>2$ доказательство аналогично (используем   индукцию).

\vspace{2mm}
\bf Относительное кратное перемешивание. \rm 
Множества вида $X\times B$ будем называть горизонтальными, а множества вида $A\times Y$  вертикальными,  подразумевая, что они измеримы. 

%\newpage
Говорим, что косое произведение $R=(S,R_x)$, $R:(X\times Y)\to(X\times Y)$, обладает относительным перемешиванием кратности $n$,
 если   для любых  горизонтальных множеств $H_0,H_1,\dots,H_n$ при $k_1,\dots, k_n\to+\infty$ для  
для индикаторов $\chi_{k_1,\dots,k_n}$ множеств 
$${H_0\cap R^{k_1}H_1\cap R^{k_1+k_2}H_2\dots\cap  R^{k_1+\dots+k_n}H_n}$$  
  выполнено 
$$ \pi \chi_{k_1,\dots,k_n}\to_\mu  const= \mu(H_0 )\mu(H_1 )  \dots \mu(H_n),$$
где $\pi$ -- ортопроекция на $L_2(X,\mu)\otimes {\bf 1}$.

\vspace{2mm}
 \bf Теорема 5.2. \it  Пусть расширение $R=(S,R_x)$, обладающее $(I,\Theta)$-коциклом,  
является  относительно  перемешивающим  с кратностью $m-1$, а автоморфизм $S$  
обладает перемешиванием кратности $m\geq 2$.   Тогда $R$ перемешивает с кратностью $m$.  
\rm

\vspace{2mm}
 \bf Лемма 5.3. \it Пусть   $R=(S,R_x)$ обладает относительным перемешиванием кратности $m-1$, а $S$ перемешивает с кратностью $m\geq 2$. Тогда для вертикальных множеств $V_0, \dots, V_m$ и горизонтальных $H_1,\dots, H_m$ при  $k_1,\dots, k_m\to+\infty$ выполнено  
$$\lambda\left(V_0\cap R^{k_1}(V_1\cap H_1) \cap \dots R^{k_1+\dots +k_m}(V_m\cap H_m)\right)\to 
\lambda(V_0 )\lambda(V_1 ) \dots \lambda(V_m)\lambda(H_1 ) \dots \lambda(H_m).$$\rm

\vspace{2mm}
Утверждение леммы непосредственно  вытекает из  ее условий.

\vspace{2mm}
Доказательство теоремы. Пусть $m=2$, рассмотрим  
 $\nu$ -- самоприсоединение автоморфизма $R$ следующего статистического происхождения.  Для некоторых последовательностей $ k_1(j),k_2(j)\to+\infty$ для всех наборов  
$\lambda$-измеримых  множеств $C_0, C_1,C_2$ ($\lambda=\mu\otimes\mu$) выполнено
$$\nu(C_0\times C_1\times C_2)= \lim_{j\to\infty}\lambda\left(C_0\cap R^{k_1(j)}C_1 \cap R^{k_1(j)+k_2(j)}C_2\right).$$ 
Заметим, что, если все такие самоприсоединения оказываются равными $\mu^3$, то $R$ обладает перемешиванием кратности 2 (это простое и давнее наблюдение позволило в ряде интересных случаев
кропотливое изучение асимптотического поведения множеств вида $C_0\cap R^{k_1(j)}C_1 \cap R^{k_1(j)+k_2(j)}C_2$ заменить на рассуждения о джойнингах).

Косое произведение $R$ удовлетворяет условиям леммы 5.3 при $m=2$, поэтому  проекция $\nu$ на $X^3$ совпадает с $\mu^3$,  а для  системы  $\nu_{(x_0,x_1,x_2)}$ условных мер на $Y^3$
имеем
$$\nu ((A_0\times B_0)\times (A_1\times B_1)\times (A_2\times B_2)) =
\int_{A_0\times A_1\times A_2} \nu_{(x_0,x_1,x_2)}(B_0\times B_1\times B_2)\,d\mu^3. $$

При $B_0=Y$ для всевозможных наборов измеримых множеств $A_0,A_1,A_2,B_1,B_2\in \B$
 с учетом относительного перемешивания (кратности 1) для $R$, применив лемму 5.3 для множеств 
$V_i=A_i\times Y$ и $H_i=X\times B_i$, получим
$$
\int_{A_0\times A_1\times A_2} \nu_{(x_0,x_1,x_2)}(Y\times B_1\times B_2)\,d\mu^3=
\mu(A_0)\mu(A_1)\mu(A_2) \mu(B_1)\mu(B_2). $$
Это влечет за собой принадлежность (почти всех) условных мер $\nu_{x_0,x_1,x_2}$  классу $M_3$. 
Далее  повторяем рассуждения из доказательства теоремы 3.2 с той лишь разницей, что вместо класса $M_4$ теперь рассматривается класс $M_3$. В результате получим
$\nu_{x_0,x_1,x_2}\equiv \mu^3$. Таким образом,
$\nu(C_0\times C_1\times C_2)=\lambda(C_0)\lambda(C_1)\lambda(C_2).$

Из сказанного вытекает, что при $k_1,k_2\to+\infty$ выполнено 
 $$\lambda\left(C_0\cap R^{k_1}C_1 \cap R^{k_1+k_2}C_2\right)\to \lambda(C_0)\lambda(C_1)\lambda(C_2),$$
т.е. $R$ обладает перемешиванием кратности 2. Для $m>1$ теорема доказывается аналогично.

\vspace{3mm}  
\bf Благoдарности. \rm  Автор весьма признателен рецензенту за  замечания и указание на пробел в первоначальном  доказательстве теоремы 4.1.

%\newpage


\normalsize
\begin{thebibliography}{99}

\bibitem{GTW}    E.~ Glasner,  J.-P.~Thouvenot,    B.~ Weiss, On some generic classes of ergodic measure preserving transformations, Тр. ММО, 82, no. 1, МЦНМО, М., 2021, 19-44; Trans. Moscow Math. Soc., 82 (2021), 15-36 

\bibitem{22} V.V.~Ryzhikov. Generic extensions of ergodic actions. arXiv:2209.09160

\bibitem{AGTW}  T.~Austin,  E.~Glasner, J.-P.~Thouvenot,   B.~Weiss,,  An ergodic system is dominant exactly when it has positive entropy,  Erg. Th. Dynam. Systems, to appear 

\bibitem{JR}  A. del Junco,  D.~Rudolph, On ergodic action whose self-joinings are graphs, 
Erg. Th. Dynam. Systems, 7 (1987), 531-557 

\bibitem{R}В\,.В.~Рыжиков, Сплетения тензорных произведений и стохастический централизатор динамических систем, Матем. сб., 188:2 (1997), 67-94; Intertwinings of tensor products, and the stochastic centralizer of dynamical systems, Sb. Math., 188:2 (1997), 237-263


\bibitem{FK} B.~Fayad,  A.~Kanigowski, Multiple mixing for a class of conservative surface
flows, Invent. math.  203 (2016), 555-614


\bibitem{RT}В\,.В.~Рыжиков,  Ж.-П.~Тувено, Дизъюнктность, делимость и квазипростота сохраняющих меру действий, Функц. анализ и его прил., 40:3 (2006), 85-89; Disjointness, divisibility, and quasi-simplicity of measure-preserving actions, Funct. Anal. Appl., 40:3 (2006), 237-240

\bibitem{R97}В.\, В.~Рыжиков, Полиморфизмы, джойнинги и тензорная простота динамических систем, Функц. анализ и его прил., 31:2 (1997),  45-57
Polymorphisms, joinings, and the tensor simplicity of dynamical 
systems, Funct. Anal. Appl., 31:2 (1997), 109-118

\bibitem{LP}  M.~Lemanczyk, F.~Parreau, Rokhlin extensions and lifting disjointness.
Ergod. Th. and Dynam. Syst.,  23:5 (2003), 1525-1550


\bibitem{L}М.\,Е.~ Липатов, Классификация коциклов над эргодическими автоморфизмами со значениями в группе Лоренца. Рекуррентность коциклов, Матем. заметки, 93:6 (2013),  869-877 
Classification of Cocycles over Ergodic Automorphisms with values in the Lorentz group and Recurrence of  Cocycles. Math. Notes, 93:6 (2013), 850-857


\bibitem{K}J.~King, Ergodic properties where order 4 implies infinite order, Israel J. Math., 80 (1992), 65-86 


\bibitem{R93}В\,.В.~Рыжиков, Джойнинги, сплетения, факторы и перемешивающие свойства динамических систем, Изв. РАН. Сер. матем., 57:1 (1993), 102-128
 Joinings, intertwining operators, factors, and mixing properties of dynamical systems, Izv. Math., 42:1 (1994), 91-114


\bibitem{T}  С.\,В.~Тихонов, О нарушении кратного перемешивания, близком к экстремальному, Тр. ММО, 82, № 1, МЦНМО, М., 2021, 205-215
 A violation of multiple mixing close to an extremal, Trans. Moscow Math. Soc., 82 (2021), 173-181


\bibitem{R96}В\,.В.~Рыжиков, Четная и нечетная простота динамических систем с инвариантной мерой, Матем. заметки, 60:3 (1996), 470-473; Even and odd primality  of dynamical systems with invariant measure, 
Math. Notes, 60:3 (1996), 353-356



\bibitem{GW}   E.~ Glasner,   B.~ Weiss, Relative weak mixing is generic, Sci. China Math. 62 (2019), 
no. 1, 69-72.


\bibitem{V} V.~Vershik, Polymorphisms, Markov processes, and quasi-similarity, Discrete Contin. Dyn. Syst., 13:5 (2005), 1305-1324



\end{thebibliography}
\end{document}